\documentclass[12pt]{article}
\usepackage{srcltx}

\linethickness{0.4pt}

\newcommand{\bi}{\bibitem}
\newcommand{\nb}{\newblock}

\newcommand{\be}[1]{\begin{equation}\label{#1}}
\newcommand{\ee}{\end{equation}}

\newcommand{\la}{\langle\,}
\newcommand{\ra}{\,\rangle}

\newcommand{\prf}{{\bf Proof.}\ }

\newcommand{\ccc}{{\cal C}}

\newcommand{\bdelta}{\bar\delta}
\newcommand{\distan}{\mathop{\rm dist}}

\newcommand{\hgt}{\mathop{\rm ht}}

\newtheorem{thm}{\quad Theorem}
\newtheorem{lm}{\quad Lemma}

\title{On the density of Cayley graphs of R.\,Thompson's group $F$ in symmetric generators}
\author{\vspace{2ex}
V. S. Guba\thanks{This work is partially supported by the Russian Foundation
for Basic Research, project no. 19-01-00591 A.}\\
Vologda State University,\\
15 Lenin Street,\\
Vologda\\
Russia\\
160600\\
E-mail: guba{@}uni-vologda.ac.ru}
\date{}

\begin{document}

\maketitle

\begin{abstract}

By the density of a finite graph we mean its average vertex degree. For an $m$-generated
group, the density of its Cayley graph in a given set of generators, is the supremum
of densities taken over all its finite subgraphs. It is known that a group with $m$ generators
is amenable iff the density of the corresponding Cayley graph equals $2m$.

A famous problem on the amenability of R.\,Thompson's group $F$ is still open. What is known
due to the result by Belk and Brown, is that the density of its Cayley graph in the standard
set of group generators $\{x_0,x_1\}$, is at least $3.5$. This estimate has not been exceeded
so far.

For the set of symmetric generators $S=\{x_1,\bar{x}_1\}$, where $\bar{x}_1=x_1x_0^{-1}$, the same
example gave the estimate only $3$. There was a conjecture that for this generating set the
equality holds. If so, $F$ would be non-amenable, and the symmetric generating set had doubling
property. This means that for any finite set $X\subset F$, the inequality $|S^{\pm1}X|\ge2|X|$ holds.

In this paper we disprove this conjecture showing that the density of the Cayley graph of $F$ in
symmetric generators $S$ strictly exceeds $3$. Moreover, we show that even larger generating set
$S_0=\{x_0,x_1,\bar{x}_1\}$ does not have doubling property.

\end{abstract}

\section*{Introduction}

Some introductory information here repeats the one of \cite{Gu04}.

The Richard Thompson group $F$ can be defined by the following infinite
group presentation

\be{xinf}
\la x_0,x_1,x_2,\ldots\mid x_j{x_i}=x_ix_{j+1}\ (i<j)\,\ra.
\ee
This group was found by Richard J. Thompson in the 60s. We refer to the
survey \cite{CFP} for details. (See also \cite{BS,Bro,BG}.) It is easy
to see that for any $n\ge2$, one has $x_n=x_0^{-(n-1)}x_1x_0^{n-1}$ so
the group is generated by $x_0$, $x_1$. It can be given by the following
presentation with two defining relations

\be{x0-1}
\la x_0,x_1\mid x_1^{x_0^2}=x_1^{x_0x_1},x_1^{x_0^3}=x_1^{x_0^2x_1}\ra,
\ee
where $a^b=b^{-1}ab$ by definition. Also we define a commutator $[a,b]=a^{-1}a^b=a^{-1}b^{-1}ab]$
and notation $a\leftrightarrow b$ whenever $a$ commutes with $b$, that is, $ab=ba$.

Each element of $F$ can be uniquely represented by a {\em normal
form\/}, that is, an expression of the form
\be{nf}
x_{i_1}x_{i_2}\cdots x_{i_s}x_{j_t}^{-1}\cdots x_{j_2}^{-1}x_{j_1}^{-1},
\ee
where $s,t\ge0$, $0\le i_1\le i_2\le\cdots\le i_s$, $0\le j_1\le j_2
\le\cdots\le j_t$ and the following is true: if (\ref{nf}) contains
both $x_i$ and $x_i^{-1}$ for some $i\ge0$, then it also contains
$x_{i+1}$ or $x_{i+1}^{-1}$ (in particular, $i_s\ne j_t$).

An equivalent definition of $F$ can be given in the following way. Let us
consider all strictly increasing continuous piecewise-linear functions
from the closed unit interval onto itself. Take only those of them that
are differentiable except at finitely many dyadic rational numbers and
such that all slopes (derivatives) are integer powers of $2$. These
functions form a group under composition. This group is isomorphic to $F$.
Another useful representation of $F$ by piecewise-linear functions can be
obtained if we replace $[0,1]$ by $[0,\infty)$ in the previous definition
and impose the restriction that near infinity all functions have the form
$t\mapsto t+c$, where $c$ is an integer.

The group $F$ has no free subgroups of rank $>1$. It is known that $F$ is
not elementary amenable (EA). However, the famous problem about amenability
of $F$ is still open. If $F$ is amenable, then it is an example of a
finitely presented amenable group, which is not EA. If it is not
amenable, then this gives an example of a finitely presented group,
which is not amenable and has no free subgroups of rank $>1$.  Note that
the first example of a non-amenable group without free non-abelian
subgroups has been constructed by Ol'shanskii \cite{Olsh}. (The
question about such groups was formulated in \cite{Day}, it is also
often attributed to von Neumann \cite{vNeu}.) Adian \cite{Ad83}
proved that free Burnside groups with $m>1$ generators of odd exponent
$n\ge665$ are not amenable. The first example of a finitely presented
non-amenable group without free non-abelian subgroups has been recently
constructed by Ol'shanskii and Sapir \cite{OlSa}. Grigorchuk \cite{Gri}
constructed the first example of a finitely presented amenable group
not in EA.
\vspace{2ex}

It is not hard to see that $F$ has an automorphism given by $x_0\mapsto x_0^{-1}$,
$x_1\mapsto x_1x_0^{-1}$. To check that, one needs to show that both defining
relators of $F$ in (\ref{x0-1}) map to the identity. This is an easy calculations using
normal forms. After that, we have an endomorphism of $F$. Aplying it once more,
we have the identity map. So this is an automorphism of order $2$.

Notice that $F$ has no non-Abelian homomorphic images \cite{CFP}. So in order to
check that an endomorphism of $F$ is a monomorphism, it suffices to show that the
image of the commutator $[x_0,x_1]=x_0^{-1}x_1^{-1}x_0x_1=x_2^{-1}x_1=x_1x_3^{-1}$
is nontrivial.

Later we will add more arguments to the importance of the symmetric set $S=\{x_1,\bar{x}_1=x_1x_0^{-1}\}$.
Obviously, it also generates $F$. It is easy to apply Tietze transormation to get a
presentation of $F$ in the new generating set from (\ref{x0-1}). So we let $\alpha=x_1^{-1}$,
$\beta=\bar{x}_1^{-1}=x_0x_1^{-1}$. It follows that $x_0=\beta\alpha^{-1}$. The first defining relation
of (\ref{x0-1}) says that $x_1^{x_0}\leftrightarrow x_1x_0^{-1}$ so
$\alpha^{\beta^{\alpha^{-1}}}\leftrightarrow\beta$. Therefore, $\alpha^{\beta}\leftrightarrow\beta^{\alpha}$.
From this relation we can derive $x_1^{x_0^2}=x_1^{x_0x_1}$ in the opposite direction.

Now the second defining relation of (\ref{x0-1}) means that $x_1^{x_0^2}\leftrightarrow x_1x_0^{-1}$, that is.
$\alpha^{\beta\alpha^{-1}\beta\alpha^{-1}}\leftrightarrow\beta$. Conjugating by $\alpha$, we get
$\alpha^{\beta\alpha^{-1}\beta}\leftrightarrow\beta^{\alpha}$. Conjugation by $\alpha$ once more implies
that $\alpha^{\beta\alpha^{-1}\beta\alpha}\leftrightarrow\beta^{\alpha^2}$. Since $\alpha^{\beta}$ commutes
with $\beta^{\alpha}=\alpha^{-1}\beta\alpha$, we conclude that the left-hand side is $\alpha^{\beta}$ so we
get the relation $\alpha^{\beta}\leftrightarrow\beta^{\alpha^2}$. Clearly, from this relation we can derive
$x_1^{x_0^3}=x_1^{x_0^2x_1}$. Therefore, by standard Tietze transformations we obtain the following presentation
of $F$ in terms of symmetric generating set:

\be{ab}
\la\alpha,\beta\mid\alpha^{\beta}\leftrightarrow\beta^{\alpha},\alpha^{\beta}\leftrightarrow\beta^{\alpha^2}\ra.
\ee

Of course, from the symmetry reasons we know that $\beta^{\alpha}\leftrightarrow\alpha^{\beta^2}$
also holds in $F$. Therefore, it is a consequence of the two relations of (\ref{ab}). Moreover, one can check that
for any positive integers $m$, $n$ it holds $\alpha^{\beta^m}\leftrightarrow\beta^{\alpha^n}$ as a consequence
of the defining relations.

\section{Density}
\label{dens}

By the {\em density\/} of a finite graph $\Gamma$ we mean
the average value of the degree of a vertex in $\Gamma$. More precisely,
let $v_1$, \dots, $v_k$ be all vertices of $\Gamma$. Let $\deg_\Gamma(v)$
denote the degree of a vertex $v$ in the graph $\Gamma$, that is, the number
of oriented edges of $\Gamma$ that come out of $v$. Then

\be{dgrform}
\delta(\Gamma)=\frac{\deg_\Gamma(v_1)+\cdots+\deg_\Gamma(v_k)}k
\ee
is the density of $\Gamma$.

Let $G$ be a group generated by a finite set $A$. Let $C(G,A)$ be the
corresponding (right) Cayley graph. Recall that the set of vertices of
this graph is $G$ and the set of edges is $G\times A^{\pm1}$. For an
edge $e=(g,a)$, its initial vertex is $g$, its terminal vertex is $ga$,
and the inverse edge is $e^{-1}=(ga,a^{-1})$. The {\em label\/} of $e$
equals $a$ by definition. For the Cayley graph $C=C(G,A)$ we define
the number

\be{denscayley}
\bdelta(C)=\sup\limits_\Gamma\delta(\Gamma),
\ee
where $\Gamma$ runs over all finite subgraphs of $C=C(G,A)$. So this number
is the least upper bound of densities of all finite subgraphs of $C$.
If $C$ is finite, then it is obvious that $\delta(C)=\bdelta(C)$. So we
may call $\bdelta(C)$ the {\em density of the Cayley graph\/} $C$.
\vspace{1ex}

This concept was used in \cite{AGL08} to study densities of the Cayley graphs of $F$.
\vspace{1ex}

Recall that a group $G$ is called {\em amenable\/} whenever there exists
a finitely additive normalized invariant mean on $G$, that is, a mapping
$\mu\colon{\cal P}(G)\to[0,1]$ such that $\mu(A\cup B)=\mu(A)+\mu(B)$
for any disjoint subsets $A,B\subseteq G$, $\mu(G)=1$, and $\mu(Ag)=
\mu(gA)=\mu(A)$ for any $A\subseteq G$, $g\in G$. One gets an equivalent
definition of amenability if only one-sided invariance of the mean is
assumed, say, the condition $\mu(Ag)=\mu(A)$ ($A\subseteq G$, $g\in G)$.
The proof can be found in \cite{GrL}.

The class of amenable groups includes all finite groups and all abelian groups.
It is invariant under taking subgroups, quotient groups, group extensions, and
ascending unions of groups. The closure of the class of finite and abelian
groups under these operations is the class EA of {\em elementary amenable\/}
groups. A free group of rank $>1$ is not amenable. There are many useful
criteria for (non)amenability \cite{Fol,Kest,Gri80}. We need to mention the
two properties of a finitely generated group $G$ that are equivalent to
non-amenability.
\vspace{1ex}

{\bf NA$_1$.}\ {\sl If $G$ is generated by $m$ elements and $C$ is the
corresponding Cayley graph, then the density of $C$ does not have the
maximum value, that is, $\bdelta(C)<2m$.}
\vspace{1ex}

Note that if NA$_1$ holds for at least one finite generating set,
then the group is not amenable and so the same property holds for
any finite generating set. For the proof of this property, we need
to use the well-known {\em F\o{}lner condition\/} \cite{Fol}. For
our reasons it is convenient to formulate this condition as follows.

Let $C$ be the Cayley graph of a group. By $\distan(u,v)$ we denote the
distance between two vertices in $C$, that is, the length of a shortest
path in $C$ that connects vertices $u$, $v$. For any vertex $v$ and a
number $r$ let $B_r(v)$ denote the ball of radius $r$ around $v$, that is,
the set of all vertices in $C$ at distance $\le r$ from $v$. For any set
$Y$ of vertices, by $B_r(Y)$ we denote the $r$-neighbourhood of $Y$,
that is, the union of all balls $B_r(v)$, where $v$ runs over $Y$. By
$\partial Y$ we denote the (outer) {\em boundary\/} of $Y$, that is, the set
$B_1(Y)\setminus Y$. The F\o{}lner condition (for the case of a finitely
generated group) says that $G$ is amenable whenever $\inf\#\partial Y/\#Y=0$,
where the infimum is taken over all non-empty finite subsets of $G$ for
a Cayley graph of $G$ in finite number of generators (this property does not
depend on the choice of a finite generating set). Any finite set $Y$ of
vertices in $C$ defines a finite subgraph (also denoted by $Y$). The degree
of any vertex $v$ in $C$ equals $2m$, where $m$ is the number of generators.
We know that exactly $\deg_Y(v)$ of the $2m$ edges that come out of $v$,
connect the vertex $v$ to a vertex from $Y$. The other $2m-\deg_Y(v)$ edges
connect $v$ to a vertex from $\partial Y$. Note that each vertex of
$\partial Y$ is connected by an edge to at least one vertex in $Y$. This
implies that the cardinality of $\partial Y$ does not exceed the sum
$\sum(2m-\deg_Y(v))$ over all vertices of $Y$. Dividing by $\#Y$ (the number
of vertices in $Y$) implies the inequality $\#\partial Y/\#Y\le2m-\delta(Y)$.
If $\bdelta(C)=2m$, then $Y$ can be chosen such that $\delta(Y)$ is
arbitrarily close to $2m$ so $\#\partial Y/\#Y$ will be arbitrarily close to
$0$. On the other hand, for any vertex $v$ in $Y$ there are at most $2m$
edges that connect $v$ to a vertex in $Y$. Therefore, the sum
$\sum(2m-\deg_Y(v))$ does not exceed $2m\#\partial Y$. So
$2m-\delta(Y)\le2m\#\partial Y/\#Y$. If the right hand side can be made
arbitrarily close to $0$, then $\delta(Y)$ approaches $2m$ so $\bdelta(C)=2m$.
\vspace{1ex}

{\bf NA$_2$.}\ {\sl If $C$ is the Cayley graph of $G$ in a finite set of
generators, then there exists a function $\phi\colon G\to G$ such that
a$)$ for all $g\in G$ the distance $\distan(g,\phi(g))$ is bounded from above
by a constant $K>0$, b$)$ any element $g\in G$ has at least two preimages
under $\phi$.}
\vspace{1ex}

An elegant proof of this criterion based on the Hall -- Rado theorem can be
found in \cite{CGH}, see also \cite{DeSS}. Note that this property also does
not depend on the choice of a finite generating set. A function $\phi$ from
NA$_2$ will be called a {\em doubling function\/} on the Cayley graph $C$.

We need a definition. Suppose that NA$_2$ holds for the Cayley graph
of a group $G$ for the case $K=1$. Then we say that the Cayley graph
$C$ is {\em strongly non-amenable\/}. The function $\phi\colon G\to G$ will
be called a {\em strong doubling function\/} on the Cayley graph $C$.
Note that each vertex is either invariant under $\phi$ or it maps into
a neighbour vertex. We know that NA$_2$ holds if and only if the group
is not amenable, that is, $\bdelta(C)<2m$. Now we would like to find
out what happens if the Cayley graph of a $2$-generated group is strongly
non-amenable.

The following fact was proved in \cite{Gu04}.
\vspace{2ex}

{\bf Theorem.} {\em
The Cayley graph of a group with two generators is strongly non-amenable
if and only if the density of this graph does not exceed $3$.}
\vspace{2ex}

It is also convenient to use the concept of {\em Cheeger boundary\/} $\partial^{\ast}Y$ of
a finite subgraph in the Cayley graph of a group regarded as a set of vertices, as above. It consists
of all directed edges that start at a vertex in $Y$ and end at a vertex outside $Y$. Clearly,
the density of $Y$ as a subgraph equals $2m\#Y-\#\partial^{\ast}Y$.

We have to mention that the density of a Cayley graph of a group is closely related to an {\em isoperimetric constant\/}
$\iota_*$ of a graph defined as $\#\partial^{\ast}Y/\#Y$; see also \cite{CGH}). Namely, one has the equality $\iota_*(C)+\bdelta(C)=2m$
for the Cayley graph $C$ of an $m$-generated group.
\vspace{1ex}

The above Theorem applied to the Cayley graph $\ccc$ of $F$ in any two generators ($x_0$, $x_1$, or $\alpha$, $\beta$) means that if we cannot
find a subgraph in with density greater than $3$, then there exists a strong doubling function on $\ccc$. One can imagine this
doubling function in the following way. Suppose that a bug lives in each vertex of $\ccc$. We allow these bugs to jump at the same time
such that each bug either returns to its initial position or it jumps to a neighbour vertex. As a result, we must have at least
two bugs in each vertex.

It is natural to ask how much the value of $\delta(Y)$ can be for the finite subgraphs we are able to construct.
In \cite{Gu04} it was constructed a family of finite subgraphs with density approaching $3$. In the Addendum yo
the same paper, there was a modification of the above construction showing that there are subgraphs with density
strictly greater than $3$. A much stronger result was obtained in \cite{BB05}. This was a family of finite
subgraphs with density approaching $3.5$. We will describe this example in the next Section. Before that, we
present a technical lemma.

First of all, we regard finite subgraphs in Cayley graphs of groups as automata, that is, labelled oriented
graphs. Let $v$ be a vertex and let $a$ be a group generator or its inverse. We say that the automaton {\em
accepts\/} $a$ whenever it has an edge labelled by $a$ starting at $v$. If the automaton does not accept $a$,
then the edge labelled by $a$ starting at $v$ in the Cayley graph, belongs to the Cheeger boundary. We claim
that the number of such edges labelled by $a$ is the same that the number of edges labelled by $a^{-1}$.

\begin{lm}
\label{invlet}
Let $G$ be a finitely generated group and let ${\cal C}={\cal C}(G,A)$ be its Cayley graph. Let $Y$ be a nonempty
finite subgraph of $\cal C$. Then for any $a\in A^{\pm1}$ the number of edges in the Cheeger boundary $\partial^{\ast}
Y$ labelled by $a$ is the same as the number of edges in $\partial^{\ast}Y$ labelled by $a^{-1}$.
\end{lm}

\prf We establish a natural bijection between edges of both types. Let $e$ be an edge labelled by $a$ in $\partial
Y$. Its starting vertex $v$ belongs to $Y$. Let $v_0=v$, and for any $n\ge0$ let $v_{n+1}$ be the starting point
of an edge in ${\cal C}$ labelled by $a$ whose terminal point is $v_n$. If $a$ has an infinite order in $G$,
then all vertices of the form $v_n$ ($n\ge0$) differ from each other. In this case, since $Y$ is finite, there
is the smallest $n > 0$ such that $v_n$ does not belong to $Y$. So $v_{n-1}$ belongs to $Y$, and the automaton
dots not accept the egde from $v_{n-1}$ to $v_n$ with label $a^{-1}$. This edge $f$ will correspond to $e$.

Suppose that $a$ has finite order in $G$. Then there is a loop in $\cal C$ at $v$ labelled by a power of $a$.
This loop has vertices outside $Y$. So, as in the previous paragraph, we can choose the smallest $n$ with the
same property. In this case we also let $e\mapsto f$, as above.

It is clear that we have a bijection between edges in $\partial^{\ast}Y$ labelled by $a$ and $a^{-1}$. The
inverse mapping $f\mapsto e$ is the same as above if we replace $a$ in the beginning by $a^{-1}$.

The proof is complete.
\vspace{1ex}

To find the density of a subgraph, we will need to know the number of edges in its Cheeger boundary. If we found
this number for a generator $a$, then we automatically know the number of edges for $a^{-1}$ due to the above
Lemma.

\section{The Brown -- Belk Construction}
\label{bbc}

Let us recall the concept of a rooted binary tree. Formally, the definition of a
rooted binary tree can be done be induction.

1) A dot $.$ is a rooted binary tree.

2) If $T_1$, $T_2$ are rooted binary trees, then $(T_1\hat{\ \ }T_2)$ is a rooted binary tree.

3) All rooted binary trees are constructed by the above rules.
\vspace{1ex}

Instead of formal expressions, we will use their formal realizations. A dot will be regarded as
a point. It coincides with the root of that tree. If $T=(T_1\hat{\ \ }T_2)$, then we draw a {\em caret\/} for $\hat{}$ as a union of two
closed intervals $AB$ (goes left down) and $AC$ (goes right down). The point $A$ is the roof of $T$. After that, we draw trees for $T_1$, $T_2$ and
attach their roots to $B$, $C$ respectively in such a way that they have no intersections. It is standard that
for any $n\ge0$ the number of rooted binary trees with $n$ carets is equal ti the $n$th Catalan number $c_n=\frac{(2n)!}{n!(n+1)!}$.

Each rooted binary trees has {\em leaves\/}. Formally they are defined as follows: for the one-vertex tree
(which is called {\em trivial\/}) the only leaf coincides with the root. In case $T=(T_1\hat{\ \ }T_2)$, the set of leaves
equals the union of the sets of leaves for $T_1$ and $T_2$. In this case the leaves are exactly vertices of degree
$1$.

We also need the concept of a {\em height\/} of a rooted binary tree. For the trivial tree, its height equals
$0$. For $T=(T_1\hat{\ \ }T_2)$, its height is $\hgt T=\max(\hgt T_1,\hgt T_2)+1$.

Now we define a {\em rooted binary forest\/} as a finite sequence of rooted binary trees $T_1$, ... , $T_m$,
where $m\ge1$. The leaves of it are the leaves of the trees. It is standard from combinatorics that the number
of rooted binary forests with $n$ leaves also equals $c_n$. The trees are enumerated from left to right and they
are drawn in the same way.

A {\em marked\/} (rooted binary) forest if the the above forest where one of the trees is marked.

Let $n\ge1$, $k\ge0$ be integer parameters. By $BB(n,k)$ we denote the set of marked forests that have $n$ leaves, and
each tree has height at most $k$. The group $F$ has a {\em left\/} partial action on this set. Namely, $x_0$ acts by
shifting the marker left if this is possible. The action of $x_1$ is as follows. If the marked tree is trivial,
this is not applied. If the marked tree is $T=(T_1\hat{\ \ }T_2)$, then we remove its caret and mark the tree
$T_1$. It is easy to see that applying $\bar{x}_1=x_1x_0^{-1}$ means the same replacing $T_1$ by $T_2$ for the marked tree.

The action of $x_1^{-1}$ and $\bar{x}_1^{-1}$ are defined analogously. Namely, if the marked tree of a forest is
rightmost, then $x_1^{-1}$ cannot be applied. Otherwise, if the marked tree $T$ has a tree $T''$ to the right of
it, then we add a caret to these trees and the tree $T\hat{\ \ }T''$ will be marked in the result. Notice that
if we are inside $B(n,k)$, then both trees $T$, $T''$ must have height $< k$: otherwise $x_1^{-1}$ cannot be
applied. For the action of $\bar{x}_1^{-1}$, it cannot be applied if $T$ is leftmost. Otherwise the marked tree $T$ has
a tree $T'$ to the left of it. Here we add a caret to these trees and the tree $T'\hat{\ \ }T$ will be marked in the
result. As above, both trees $T'$, $T'$ must have height $< k$ to be possible to stay inside $B(n,k)$.

We have to emphasize that the definition of these actions is very important for Section~\ref{mr}. So the reader
has to keep in mind these rules. We will use them without reference.

It can be checked directly that applying defining relations of $F$ leads to the trivial action (in case when
the action of each letter is possible. For details we refer to \cite{BB05}. So one can regard $BB(n,k)$ as a
set of vertices of the Cayley graph of $F$. This can be done for each of the three generating sets
$\{x_0,x_1\}$, $\{x_1,\bar{x}_1\}$, and $\{x_1,\bar{x}_1,x_0\}$.

For any fixed $k$, let $n\gg k$. Since any tree of height $k$ has at most $2^k$ leaves, any forest in $B(n,k)$
contains at least $\frac{n}{2^k}$ trees. Therefore. if we randomly take a marked forest, the probabililty for
this vertex of an automaton to accept both $x_0$, $x_0^{-1}$ approaches $1$. Now look at the probability to
accept $x_1^{-1}$. The contrary holds if and only if the marked tree is trivial. We may assume this tree is
not the rightmost one of the forest. Then we remove the trivial tree and move the marker to the right. As a
result, we obtain an element of $B(n-1,k)$. The inverse operation is always possible. So the probability we
are interested in, equals $\#B(n-1,k)/\#B(n,k)$. It approaches some number $\xi_k$ as $n\to\infty$. If $k$ is
big enough, then $\xi_k$ is close to $\frac14$. Indeed, for large $k$ the number of elements in the set $B(n,k)$
grows almost like $4^n$, as Catalan numbers do.

For the inverse letter $x_1^{-1}$, straightforward estimating the probability not to accept it is more complicated.
However, it is the same as for $x_1$ due to Lemma \ref{invlet}. We see that the number of outer edges in the
subgraph (that is, the edges in its Cheeger boundary of $B(n,k)$ approaches one half of the cardinality of this
set. This means that the density of the set $B(n,k)$ approaches $3.5$.

To be more precise, let us add some calculations. First of all, let $\Phi_k(z)$ be the generating function of the
set of rooted binary trees of height $le k$ with $n$ leaves. Clearly, $\Phi_0(z)=z$. For $k > 0$ we have either
a trivial tree that correspond to the summand $z$, or it has an upper caret. Removing it, we have an ordered
pair of trees of height $\le k-1$. Hence $\Phi_k(z)=z+\Phi_{k-1}(z)^2$.

So we have a sequence of polynomials with positive integer coefficients. All of them are increasing functions on
$z\ge0$ and approach infinity as $z\to\infty$. So there exists a unique solution of the equation $\Phi_k(z)=1$.
We denote it by $\xi_k$. This is a decreasing sequence. Let us show that $\xi_k\to\frac14$ as $k\to\infty$.

First of all, by induction on $k$ one can easily check that $\Phi_k(\frac14) < \frac12$ for all $k\ge0$. Thus
$\frac14 < \xi_k$. On the other hand, every tree with $n\le k$ carets (so $n+1$ leaves) has height $\le k$.
Hence the first terms of $\Phi_k(z)$ coincide with Catalan numbers: the coefficient on $x_{n+1}$ equals $c_n$
for $n\le k$. It is known that the series $\Phi(z)=c_0z+c_1z^2+\cdots+c_nz^{n+1}+\cdots=\frac{1-\sqrt{1-4z}}2$
has radius of convergence $\frac14$. So for any $z > \frac14$, the partial sums of the series approach infinity.
Thus $c_0z+c_1z^2+\cdots+c_nz^{n+1} > 1$ if $n$ is sufficiently large. In particular, $\Phi_k(z) > 1$ whenever
$k$ is large enough. So $\frac14 < \xi_k < z=\frac14+\varepsilon$ for $k\gg1$. This proves what we claim.

\section{Main Results}
\label{mr}

\begin{thm}
\label{denssym}
The density of the Cayley graph of Thompson's group $F$ in symmetric generating set $S=\{x_1,\bar{x}_1=x_1x_0^{-1}\}$
is strictly greater than $3$.
\end{thm}

\prf First we consider the Brown -- Belk set $B(n,k)$. It gives a subgraph in the Cayley graph $\mathcal C$ of
the group $F$ in generating set $S$. Let us find the generating function of
this set for any $k$. The coefficient on $z^n$ will show the number of marked forests with $n$ leaves where
all trees of this forest have height $\le k$.

The marked tree of the forest has generating function $\Phi_k(z)$. To the left of it, we may have any number
of trees including zero. Thus for this part we get generating function
$1+\Phi_k(z)+\Phi_k^2(z)+\cdots=\frac1{1-\Phi_k(z)}$. The same for the trees to the left of the marker.
Therefore, we get a function $\Psi_k(z)=\frac{\Phi_k(z)}{(1-\Phi_k(z))^2}$. Its coefficient on $z^n$ in the
series expansion is exactly the cardinality of $B(n,k)$. We shall denote it by $\beta_{nk}$.

The radius of convergence of the series for $\Psi_k(z)$ equals $\xi_k^{-1}$. On the other hand, the quotient
$\frac{\beta_{n-1,k}}{\beta_{nk}}$ approaches the reciprocal of the radius, that is, for any $k$ one has
\be{appr}
\frac{\beta_{n-1,k}}{\beta_{nk}}\to\xi_k
\ee
as $n\to\infty$.

The automaton corresponding to $B(n,k)$ does not accept $x_1$ whenever the marked tree is trivial or it is the
rightmost one in the forest. The former case happens with probability $\le\frac{2^k}n$ since each forest has
at least $\frac{n}{2^k}$ trees. So for any $k$, the probability to be rightmost is $o(1)$ as $n\to\infty$. If
the tree is not rightmost, then we remove the trivial tree and move the marker right. The number of these new
marked trees we obtain is exactly $\beta_{n-1,k}$. Indeed, an inverse operation of inserting the trivial tree
and moving the marker left is always possible. So the total probablility for a marked forest in $B(n,k)$ not
to accept $x_1$ equals $\frac{\beta_{n-1,k}}{\beta_{nk}}+o(1)=\xi_k+o(1)$ according to~(\ref{appr}).

As for symmetric generator $\bar{x}_1$, the probability has exactly the same value (replace ``rightmost" by
``leftmost"). Lemma~\ref{invlet} allows us to conclude that the same happens for inverse letters. Therefore,
the cardinality of the Cheeger boundary of $B(n,k)$ divided by the cardinality of $B(n,k)$ itself, is
$4\xi_k+o(1)$. It approaches $1$ as $k\to\infty$, so the density of $B(n,k)$ approaches $3$. At the present
time, the sets $B(n,k)$ give the best density estimate for the generating set $\{x_0,x_1\}$. So there was
a conjecture that for $S$ the density $3$ could be an optimal value.

However, this conjecture is not true. There is an essential difference between the standard generating set
and the symmetric one. If we take a random marked forest, it always accepts both $x_0$ and $x_0^{-1}$ if
the marked tree is neither leftmost nor rightmost. We already know that the probability to be leftmost
(rightmost) does not exceed $\frac{2^k}n$ so it is almost zero for $n\gg1$. Thus the degree of any vertex
of a graph is at least $2$ for almost all cases if we work with standard generating set.

Now look at the vertices of the Cayley graph $\mathcal C$. It turns out that they can be isolating. Indeed,
let we have a marked forest $\dots, T', T, T'', \dots$ where $T$ is marked. Suppose that $T$ is trivial.
Then the vertex does not accept $x_1$ as well as $\bar{x}_1$ (we cannot remove a caret). Additionally suppose
that both trees $T'$, $T''$ have height $k$. This means that we cannot apply neither $\bar{x}_1^{-1}$ (adding a
caret to $T'$ and $T$), nor $x_1^{-1}$ (adding a caret to $T$ and $T''$). So in this case we get an isolated
vertex.

What is the probability for a vertex (that is, a random marked tree from $B(n,k)$) to be isolated? If it is
small, then we have no profit from that. But it turns out that the probability is uniformly positive. That is,
there exists a global positive constant $p_0 > 0$ such that the probability of a vertex to be isolated will
be at least $p_0$ for all our graphs.

The fact we claim is sufficient to prove the theorem. Indeed, if we remove the isolated vertices from $B(n,k)$,
then we get a subgraph, say, $B'(n,k)$, where the number of its edges is the same and the number of vertices
will be less than $(1-p_0)\beta_{nk}$. Since the density is an average degree of a vertex, then the density
of the new subgraph will be at least $\frac1{1-p_0}$ multiplied by the density of $B(n,k)$, which is
$3-\varepsilon$ for arbitrarily small $\varepsilon > 0$. This means that we can approach density $\frac3{1-p_0}
> 3$ of the Cayley graph $\mathcal C$.

So let us show that the value $p_0=\frac1{260}$ can be established (so that the density of $\mathcal C$ will exceed
$3.011$). (Recall that strict inequality here makes useless the idea to find any kind of a ``doubling structure"
on $\mathcal C$, in the sense we have mentioned in the Introduction.) Direct calculations with generating
functions do not give us a clear way to prove the statement. So we will prefer a probabilistic approach.

Let $\dots, T_{-1}, T_0, T_1, \dots$ be a random marked forest. Assume that all the three trees
$T_{-1}$, $T_0$, $T_1$ are trivial. What is the probablility of that? If $T_1$ is rightmost, then we already know that the
probability is $o(1)$ so we can ignore this case. If we remove the three trivial trees and move the marker to
the tree that goes after $T_1$ (let it be $T_2$ in the above notation), then we obtain a marked forest from $B(n-3,k)$. The inverse operation
is always possible to do. So our probability is $\frac{\beta{n-3,k}}{\beta_{nk}}+o(1)=\xi_k^3+o(1)$ as $n\to\infty$.

Now we start add carets. The first one is added to $T_1$ and $T_2$. Then we add a caret to obtain $((T_1 \hat{\ \
}T_2)\hat{\ \ }T_3)$ and so on. At some step we will not be able to add a new caret. This happens if we reach
the rightmost tree (for what case the probability is very small), or we cannot add a new caret to two trees
because at least one of them has height $k$. To be more precise, let us assume that the trees $T_2$, ... ,
$T_{k+1}$ do exist in our marked forest. If not, the probability for a marked tree $T$ be close to the right
border does not exceed $\frac{(k+1)2^k}{n}=o(1)$ as $n\to\infty$. So the process of adding carets to the right
of $T$ will get us $... , T, T_1'',T_2'',...$ where at least one of the trees $T_1''$, $T_2''$ has height $k$.

The same process can be done to the left of $T$. There we get $... , T_2',T_1',T,...$, where at least one of the trees
$T_2'$, $T_1'$ has height $k$.

Suppose that both $T_1'$, $T_1''$ have height $k$. Then the marked forest $...,T_1',T,T_1'',...$ gives an
isolated vertex as we have seen before. If $T_1'$ does not have height $k$ then $T_2'$ has height $k$ so we can
swap the trees $T_2'$ and $T_1'$ in the forest. Both of these forest will have the same probability. Also if
$T_1''$ does not have height $k$ then $T_2''$ has height $k$ and we swap these trees. Then the probability of
our event (when $T_1'$ and $T_1''$ have height $k$) is at least $\frac14$ of the probablility of the event:
($T_2'$ OR $T_1'$ has height $k$) AND ($T_1''$ OR $T_2''$ has height $k$). The former is $\xi_k^3+o(1)$ since
the process of adding carets is unique and the inverse operations are possible to do. This will lead back to
the case of three trivial trees for which the probability is already known.

So we proved that the probability of a random vertex to be isolated is at least $\frac14\xi_k^3+o(1)$. It
approaches $\frac1{4^4}=\frac1{256} > \frac1{260}=p_0$. This completes the proof.
\vspace{2ex}

At the end of this Section we will obtain one more result. Let us add $x_0$ to the generating set $S$. We will
get three generators $\{x_1,\bar{x}_1,x_0\}$. What is the density of the Cayley graph here, is not known. We
only know that the isoperimetric constant $\iota^{\ast}$ is at least $1$ but we cannot prove or disprove the
strict inequality. The idea to remove isolated vertex does not work here since in the new graph the former
isolated vertex will have degree $2$ because of edges labelled by $x_0^{\pm1}$.

However, we can say something about the outer boundary $\partial$ instead of the Cheeger boundary
$\partial^{\ast}$. The question from the previous paragraph is equivalent to the following: is there a finite
set $Y\subset F$ such that $\#\partial^{\ast}Y < \#Y$? We do not know the answer but we are able to prove the
following.

\begin{thm}
\label{double}
For the symmetric generating set $S=\{x_0,x_1,\bar{x}_1=x_1x_0^{-1}\}$, there exists finite subsets $Y\subset F$ in the
Cayley graph of Thompson's group $F$ such that $\#\partial Y < \#Y$.

Equivalently, the generating set $S$ does not have doubling property, that is, there are finite subsets $Y$ in
$F$ such that the $1$-neighbourhood $\mathcal N_1(Y)=Y(\{1\}\cup S)$ has cardinality strictly less than $2\#Y$.
\end{thm}

\prf
The second statement follows to the first one (and in fact it is equivalent) since the $1$-neigbourhood of $Y$
is the disjoint union of $Y$ itself and its outer boundary.

The proof of the first statement will be easier that the proof of Theorem~\ref{denssym} since in this case it
suffices to take $Y=B(n,k)$. For every vertex $v$ in its outer boundary, we choose an edge $e$ connecting it to a
vertex $u$  in $Y$. If there are several edges with this property, we fix one of them. The aim is to estimate the number of
fixed edges, which is equal to $\partial Y$.

If the label of a fixed edge is $x_0^{-1}$, we already know that the probability is $o(1)$. Here we think in
terms of probabilities dividing the number of edges by $\#Y$.

Suppose that the edge $e$ has label $x_1^{-1}$. Then $u$ as a vertex of the automaton $Y$ does not accept $x_1$.
This means that we cannot remove a caret of the marked tree corresponding to $u$. This means that the tree is
empty. So the number of edges $e$ with label $x_1^{-1}$ does not exceed the number of marked forests with
trivial marked tree. In terms of probabilities, this gives the estimate $\xi_k+o(1)$. Exactly the same holds
for edges $e$ labelled by $\bar{x}_1^{-1}$ because of symmetry.

Now look at the number of vertices $v$ in the outer boundary for which the label of $e$ is $x_1$ or $\bar{x}_1$.
The vertex $v$ can be represented as a marked forest. After we apply $x_1^{-1}$ or $\bar{x}_1^{-1}$ to it
removing the upper caret, we get to $u$ which is a forest with all trees of height $\le k$. Therefore, the tree
$T=T_1\hat{\ \ }T_2$ which is marked for the vertex $v$, has height $k+1$. Applying $x_1^{-1}$ to it means that
the caret is removed and the marked tree becomes $T_1$.

So the vertices $v$ in the outer boundary for which the label of $e$ is $x_1$ or $\bar{x}_1$ are connected to
a vertex in $Y$ by both of these edges. So in the process of choosing edges, we may assume that the label of $e$
is $x_1$. Hence the number of vertices with this property does not exceed the number of vertices in $Y$ for
which $x_1^{-1}$ cannot be applied. By Lemma~\ref{invlet}, this number equals the one for the generator $x_1$.
The probability for that is $\xi_k+o(1)$ already know.

Summing the numbers, we obtain that $\#\partial Y/\#Y=3\xi_k+o(1) < 1$ for $k\gg1$. In fact, the constant $2$ in
the statement of the Theorem can be replaced by $\frac74+\varepsilon$ for any positive $\varepsilon$.
\vspace{1ex}

The proof is complete.
\vspace{2ex}

Notice that in order to prove non-amenability of Thompson's group $F$ (if it is the fact), it suffices to find a kind of a
doubling structure on the Cayley graph of the group. If the generating set is ``small" then we have no chances
to find this structure as the above results show. If it is very ``large", then it is more difficult to work with
the graph. So we would like to offer the generating set $\{x_0,x_1,x_2\}$ for which there are chances to find a
doubling structure on the Cayley graph.

\end{document}